\documentclass[12pt]{article}
\usepackage{amsthm, amssymb, amsmath}
\newtheorem{theorem*}{Theorem}

\def\QQ{\mathbb{Q}}
\def\ZZ{\mathbb{Z}}

\title{Unramified Alternating Extensions of Quadratic Fields}
\author{Kiran S. Kedlaya}
\date{\today}

\begin{document}

\maketitle

\begin{abstract}
We exhibit, for each $n \geq 5$, infinitely many quadratic number fields
admitting unramified degree $n$ extensions with prescribed signature
whose normal closures have Galois group $A_n$. 
This generalizes a result of Uchida and Yamamoto, which did not include
the ability to restrict the signature, and a result of Yamamura, which
was the case $n=5$.
\end{abstract}

It is a folk conjecture that for $n \geq 5$,
all but finitely many quadratic number
fields admit unramified extension fields of degree $n$ whose
normal closures have Galois group $A_n$, the alternating group on
$n$ symbols.
Uchida \cite{bib:uchida}
and Yamamoto \cite{bib:yamamoto}
proved independently that there exist infinitely
many real and infinitely many imaginary fields with unramified
$A_n$-extensions. Our main theorem is the following refinement of
this result; the case $n = 5$ was previously obtained by
Yamamura \cite{bib:yamamura} using a different argument.

\begin{theorem*}
For $r = 0, 1, \dots, \left\lfloor n/2 \right\rfloor$, there exist
infinitely many real quadratic fields admitting an unramified
degree $n$ extension with Galois group $A_n$ and having exactly
$r$ complex embeddings. In fact, the number of such real
quadratic fields with discriminant at most $N$ is at least $O(N^{1/{n-1}})$.
Moreover, these assertions remain true if we require all fields involved
to be unramified over a finite set of finite places of $\QQ$.
\end{theorem*}
\begin{proof}
The idea is to construct monic polynomials $P(x)$ with integer
coefficients and squarefree discriminant $\Delta$, so that $\QQ[x]/(P(x))$
is unramified over $\QQ(\sqrt{\Delta})$. To do so, we construct
$Q(x) = n (x-a_1)\cdots (x-a_{n-1})$ and set
$P_b(x) = b + \int_0^x P(t)\,dt$.
Then the discriminant $\Delta_b$ of $P_b(x)$ factors as
\[
\Delta_b = n^n \prod_{i=1}^{n-1} P_b(a_i) = n^n \prod_{i=1}^{n-1} (P_0(a_i) + b);
\]
since each factor is linear in $b$, a simple sieving argument
shows that the discriminant is squarefree for a positive proportion
of tuples $(a_1, \dots, a_{n-1}, b)$ in suitable ranges.

In order to make this program work, we must add constraints on the
$a_i$ and $b$ to ensure that various conditions are met. We first
make sure that $P(x)$ has integer coefficients. It suffices to
require that $a_1$ be of the form $q/n$ with $q$ an integer coprime to $n$ but
divisible by $n-1$ and by all primes less than $n$ not dividing $n$,
that $a_2, \dots, a_{n-1}$ be integers divisible by $n!$, and that $b$ be
an integer coprime to $n!$. Now
\[
Q(x) \equiv (nx-q) x^{n-2} = nx^{n-1} - q x^{n-2} \pmod{n!},
\]
and so the coefficient of $x^{m-1}$ is divisible by $m$ for $m=1, \dots, n$.
Consequently, $P_b(x)$ has integer coefficients.

Next, we force $P_b(x)$ to have a specific number of real roots,
which determines the number of real embeddings of the field $\QQ[x]/(P_b(x))$.
By homogeneity, the number of real roots of $P_b$ depends only on the
tuple $(a_1/b^{1/n}, \cdots, a_{n-1}/b^{1/n})$. Thus we can construct
intervals $(c_i, d_i)$ and $(c,d)$ such that if $a_i/m^{1/n} \in (c_i, d_i)$
for all $i$ and $b/m \in (c,d)$
for some $m$, then $P_b(x)$ has the desired number of real roots.

Next, we ensure that the splitting field of $P_b(x)$ over $\QQ$ has Galois
group $S_n$. We do this by imposing congruence conditions modulo some
auxiliary primes. Pick any degree $n$ polynomial $R(x)$
over $\ZZ$ with Galois group
$S_n$ such that the splitting fields of $R(x)$ and $R'(x)$ are linearly
disjoint. Then by the Cebotarev Density Theorem,
there exist infinitely many primes $p_1$ and $p_2$
such that $R'(x)$ factors completely
modulo $p_1$ and $p_2$, while $R(x)$ is irreducible over $p_1$ and factors
into one quadratic factor and $n-2$ linear factors over $p_2$. We may 
impose congruence conditions on the $a_i$ and $b$ so that $P_b(x)
\equiv R(x) \pmod{p_1p_2}$, which forces $P_b(x)$ to have Galois group
$S_n$ over $\QQ$.

Let us rewrite the factorization of $\Delta_b$ as
\[
\Delta_b = (n^n P_0(q/n) + bn^n) \prod_{i=2}^{n-1} (P_0(a_i) + b).
\]
The first term can written as $q^n$ plus a multiple of $n$, so is coprime
to $n$; the remaining terms are each $b$ plus a multiple of $n$, so are also
coprime to $n$. Hence $\Delta_b$ is coprime to $n$. Also, if $p$ is a prime
less than $n$ not dividing $n$, then $P_b(a_i) \equiv a_i^n + b \pmod{p}$
and so none of the factors of $\Delta_b$ is divisible by $p$ either.

For future convenience, we restrict $a_1, \dots, a_{n-1}$ to a very
special form. We require them to be of the form $A_1 \ell, \dots, A_{n-1}\ell$
for $A_1, \dots, A_{n-1}$ fixed once and for all and $\ell$ a prime.
Then by homogeneity, we can write $P_0(a_i) = B_i \ell^n$ for some integers
$B_1, \dots, B_{n-1}$. By imposing congruence conditions on $\ell$ and
$b$ modulo the primes dividing $\prod_{i<j} B_i - B_j$, we may ensure that
no prime except possibly $\ell$ divides more than one of the factors
$P_0(a_i)+b$.

Finally, we ensure that each factor of $\Delta_b$ is squarefree; this step
is analogous to the proof that $6/\pi^2$ of the positive integers are
squarefree. (We have followed \cite{bib:hooley} in this stage of
the argument.)
Fix $N$ and pick a prime $\ell$ such that $a_i/N^{1/n(n-1)}
\in (c_i, d_i)$ for all 
$i$; we will sieve over integers $b$ such that $b/N^{1/(n-1)} \in (c,d)$.
As noted above, the only prime that can divide more than one
of the factors $P_0(a_i)+b$ and $P_0(a_j)+b$ is $\ell$. Thus 
we must exclude $n(d-c)N^{1/(n-1)}/\ell + O(1)$ of the possible values $b$
in the range of interest.

Under this restriction, the factors $P_0(a_i)+b$ are pairwise coprime,
so it suffices to ensure that each one is squarefree for a positive proportion
of $b$ among the values of interest. Let $S$ denote the set of $b$
for which $b/N^{1/(n-1)} \in (c,d)$ and no $P_0(a_i)+b$ is divisible by
$\ell$. Let $N_0 = |S|$, let $N_1$ denote the number of $b \in S$ such that
each $P_0(a_i)+b$ is squarefree,
let $N_2$ denote the number of $b \in S$
such that no $P_0(a_i)+b$ is divisible by the
square of any prime less than $\xi = \frac{1}{4} \log N^{1/(n-1)}$,
and let $N_3$ denote the number of $b \in S$ such that
$P_0(a_i)+b$ is divisible by the square of a prime greater than
$\frac{1}{4} \log N^{1/(n-1)}$. These are related by the equation
\[
N_1 = N_2 + O(N_3).
\]
Now $N_2$ can be written, by inclusion-exclusion, as a sum over
squarefree numbers $l$ whose prime factors are all less than $\xi$.
Any such number is at most $N^{1/2(n-1)}$, so, with $\mu(l)$ denoting
the M\"obius function at $l$ and $d(l)$ the number of divisors of $l$,
we have
\begin{align*}
N_2 &= \sum_l n \mu(l) \left( \frac{N_0}{l} + O(1) \right)k\\
&= N_0 \sum_l \frac{n \mu(l)}{l} + O\left( \sum_{l \leq N^{1/2(n-1)}}
d(l) \right) \\
&= N_0 \prod_p \left( 1 - \frac{n}{p} \right) + O \left( \frac{x}{\log x}\right).
\end{align*}
As for $N_3$, we have the estimate
\begin{align*}
N_3 &= \sum_{\xi < p < N^{1/2(n-1)}} n \left( \frac{N_0}{p} + O(1) \right) \\
&= O\left( N_0 \sum_{\xi < p < N^{1/2(n-1)}} \frac{1}{p} \right)
+ O\left( \sum_{\xi < p < N^{1/2(n-1)}} 1 \right) \\
&= O\left(\frac{N_0}{\log N^{1/(n-1)}}\right).
\end{align*}
Putting this together, we conclude that a positive proportion of $b \in S$
yield squarefree $\Delta_b$.

We now have produced $O(N^{1/(n-1)}$ unramified $A_n$-extensions
of prescribed signature over quadratic fields of discriminant at most
$N$. Moreover, the number of distinct values of $\Delta_b$ occurring is
also $O(N^{1/(n-1)})$.
Thus at least this many quadratic fields of discriminant
less than $N$ admit unramified $A_n$-extensions of the desired signature.
\end{proof}

We have not attempted to obtain the best possible exponent in the theorem;
by varying $a_1, \dots, a_{n-1}$, one ought to be able to get an
exponent of $2/n$. It may be possible to do even better by allowing
$a_1, \dots, a_{n-1}$ to lie not in $\QQ$ but in a larger number field.

\subsection*{Acknowledgments}
Thanks to Manjul Bhargava and Noam Elkies for helpful discussions.
The author was supported by an NSF Postdoctoral Fellowship, and conducted
this research at the Mathematical Sciences Research Institute.


\begin{thebibliography}{99}

\bibitem{bib:hooley}
C. Hooley, On the power free values of polynomials, \textit{Mathematika}
\textbf{14} (1967), 21--26.

\bibitem{bib:uchida}
K. Uchida, Unramified extensions of quadratic number fields, II,
\textit{T\^ohoku Math. J.} \textbf{22} (1970), 220--224.

\bibitem{bib:yamamoto}
Y. Yamamoto, On unramified Galois extensions of quadratic number fields,
\textit{Osaka J. Math.} \textbf{7} (1970), 57--76.

\bibitem{bib:yamamura}
K. Yamamura, On unramified Galois extensions of real quadratic number
fields, \textit{Osaka J. Math.} \textbf{23} (1986), 471--478.

\end{thebibliography}
\end{document}